\RequirePackage{fix-cm}
\documentclass[smallextended]{svjour3}       % onecolumn (second format)
\smartqed  % flush right qed marks, e.g. at end of proof
\usepackage{amsmath}
\usepackage{graphicx}
\usepackage{amsfonts}
\usepackage{amssymb}
\DeclareMathOperator{\dist}{dist}
\DeclareMathOperator{\prox}{Prox}
\DeclareMathOperator{\dom}{dom}
\begin{document}

\title{Nonconvex Proximal Incremental Aggregated Gradient Method with Linear Convergence%\thanks{Grants or other notes
}


\author{Wei Peng$^1$         \and
        Hui Zhang$^1$ \and Xiaoya Zhang$^1$ %etc.
}


\institute{Wei Peng
              \\\email{weipeng0098@126.com}           %  \\
				\\\\
           Hui Zhang(Corresponding author) 
           \\\email{h.zhang1984@163.com}
           \\\\
           Xiaoya Zhang
           \\\email{zhangxiaoya09@nudt.edu.cn}
           \\\\
              $^1$ Department of Mathematics, National University of Defense Technology \\
}

\date{Received: date / Accepted: date}

\titlerunning{Nonconvex PIAG with Linear Convergence}
\authorrunning{Wei Peng et al.}
\maketitle

\begin{abstract}
In this paper, we study the proximal incremental aggregated gradient(PIAG) algorithm for minimizing the sum of L-smooth nonconvex component functions and a proper closed convex function. 
By exploiting the L-smooth property and with the help of an error bound condition, we can show that the PIAG method still enjoys some nice linear convergence properties even for nonconvex minimization. 
To illustrate this, we first demonstrate that the generated sequence globally converges to the stationary point set. 
Then, there exists a threshold such that the objective function value sequence and the iterate point sequence are R-linearly convergent when the stepsize is chosen below this threshold.

% \PACS{PACS code1 \and PACS code2 \and more}u
\end{abstract}

\section{Introduction}
\label{intro}
A fundamental optimization model emerges in numerous problems including machine learning, signal processing, image science, communication systems, and distributed optimization. Typically, the model is to minimize the sum of $N$ differentiable functions $f_i$ which are possibly nonconvex and a convex nonsmooth function $h$:
\begin{align}\label{model}
\min_{x\in\mathbb{R}^d} F(x):=\sum_{i=1}^N f_i(x)+h(x).
\end{align}
This problem often arises in large-scale, distributed, parallel optimization subfields with large $N$. 
%The popular forward-backward splitting (FBS) scheme to solve this model consists of a gradient (forward) step of $\sum_{i=1}^Nf_i(x)$ and a proximal (backward) step of $h(x)$. In this case, the gradient step of FBS to evaluate $\nabla f(x)$ requires each $\nabla f_i(x)$ to be calculated, which might be prohibitive for some large $N$. 
Directly computing $\sum_{i=1}^N\nabla f_i(x_k)$ in the popular forward-backward splitting (FBS)\cite{combettes2005signal} scheme might be prohibitive for large $N$. Thereby a natural method to approximate the gradient of $f(x)$ embedding into FBS scheme is proposed, named as the proximal incremental aggregated gradient (PIAG) method. 
The key idea of PIAG is to construct an ``inexact gradient" $g_k$ to substitute the ``exact" $\sum_{i=1}^N\nabla f(x_k)$. PIAG method is the iterative procedure of three steps:
\begin{align}
&g_k=\sum_{i=1}^N\nabla f_i(x_{k-\tau_k^i}),\\
&y_k=x_k-\alpha\cdot g_k\label{yk},\\
&x_{k+1}=\arg\min_{x\in\mathbb{R}^d}\left\{h(x)+\frac{1}{2\alpha}\|x-y_k\|^2\right\}\label{xkp1},
\end{align}
where $\tau_k^i$ are some nonnegative integers representing delayed iterations. In addition, we assume that $\tau_k^i$ never exceeds a given integer $\tau\geq0$. Thereby the exact gradient $\nabla f_i(x_k)$ is approximated by previous gradient components $\nabla f_i(x_{k-\tau_k^i})$ no more than $\tau$ iterations before. We can rewrite (\ref{yk}) and (\ref{xkp1}) into the following subproblem:
\begin{align}
x_{k+1}=\arg\min_{x\in\mathbb{R}^d}\left\{h(x)+\langle g_k,x-x_k\rangle+\frac{1}{2\alpha}\|x-x_k\|^2\right\}\label{intro_subprob}.
\end{align}
Note that under the condition of no delays involved, i.e. $\tau_{k}^i\equiv 0$, we have $g_k=\sum \nabla f_i(x_k)$, which is exactly the classic FBS scheme.

On one hand, PIAG has been investigated in several works under the convex settings. \cite{RN63} is the first to establish a global linear convergence rate of PIAG for strongly convex minimization, which guarantees that PIAG returns an $\varepsilon$-optimal solution after $\mathcal{O}(Q\tau^2\log_2(Q\tau) \log(1/\varepsilon))$ iterations, where $Q$ is the condition number. \cite{PIAG2016} showed a global linear convergence rate in $\|x_k-x^\ast\|$ with complexity no more than $\mathcal{O}(\log(1/\varepsilon)Q\tau^2)$; general distance functions are also involved in their analysis. Combining \cite{RN63} and \cite{PIAG2016}, \cite{RN64} proposed stronger linear convergence rate that achieving an $\varepsilon$-optimal solution of the function values requires at most $O(Q\tau \log(1/\varepsilon))$ iterations. \cite{zh1} gave the global linear convergence of PIAG under several strictly weaker assumptions, novel variants with better convergence rate as well as an improved rate result under strongly convex condition. \cite{zhang2017inertial} proposed an accelerated globally linearly convergent scheme under quadratic growth condition, which combines the heavy ball method with Nesterov-like acceleration.

On the other hand, there are a few studies of nonconvex FBS. A nonconvex nonsmooth version of FBS was analyzed in \cite{attouch2013convergence}, which is involved in a fundamental approach under the Kurdyka-{\L}ojasiewicz(KL) condition. \cite{bolte2014proximal,attouch2010proximal} considers a structured functions of the type $L(x,y)=f(x)+Q(x,y)+g(x)$, with the case of nonconvex FBS included. Instead of using the KL conditions, this paper is consistent with a series of studies \cite{luo1992linear,luo1993convergence,beck2006linearly,TP1,TP3} using the proximal error bound condition, under which the recent work \cite{peng} gave local linear convergence results of an accelerated nonconvex proximal gradient method. The accelerated algorithm is exploiting historical information in essence, sharing the similar viewpoint with PIAG to an extent.
  Thus inspired by their work, we analyze the linear convergence of nonconvex PIAG with the proximal error bound condition. 

\paragraph{Main contribution.} In this study, we mainly focus on the convergence analysis of PIAG for minimizing a class of nonconvex problems, under the proximal error bound condition.
First, we prove the sequence $\{x_k\}$ generated by PIAG is globally convergent to the stationary point set(Theorem 1(i)) of (\ref{model}). 
Then, when we choose the stepsize $\alpha$ below a certain threshold, the objective function value sequence is proved to be R-linearly convergent to the function value at a certain stationary point(Theorem 1(ii)). Finally, with the proved R-linearly convergent property of function value sequence, we show the iterate sequence $\{x_k\}$ generated by PIAG R-linearly converges to a certain stationary point.

The rest of the paper is structured as follows. Section 2 introduces notations and assumptions to be used. Section 3 gives the convergence analysis and section 4 concludes the paper.

%nonconvex accelerated PG results
\section{Notations \& Assumptions}
Throughout this paper, $d$-dimensional Euclidean space is denoted by $\mathbb{R}^d$ and its inner product is represented by $\langle\cdot,\cdot \rangle$. The $l_2$-norm is denoted by $\|\cdot\|$. For a nonempty closed set $\mathcal{C}\subset\mathbb{R}^d$, the distance from $x$ to $\mathcal{C}$ is represented by $\dist(x,\mathcal{C})$, where $\dist(x,\mathcal{C})=\inf_{y\in\mathcal{C}}\|x-y\|$. 
The domain of an extended-value function $h:\mathbb{R}^d\rightarrow[-\infty,+\infty]$ is defined as
$\dom h=\{x\in\mathbb{R}^d,h(x)<+\infty\}$ and $h$ is said to be proper if $h$ is never equals $-\infty$ and $\dom h\neq \emptyset$.
The gradient operator of a differentiable function is denoted by $\nabla$. The subdifferential of a proper lower closed convex function $h$ is defined as
\begin{align}
\partial h(x)=\{v\in\mathbb{R}^d:h(u)-h(x)-\langle v,u-x\rangle\geq0,~~\forall u\in\mathbb{R}^d\},
\end{align}
where $\partial h(x)$ is always a closed convex set.
The proximal operator of a proper closed function $h$ at $y\in\mathbb{R}^d$ is defined as
\begin{align}
\prox_h(y)=\arg\min_{x\in\mathbb{R}^d}\left\{h(x)+\frac{1}{2}\|x-y\|^2\right\}.
\end{align}
The sequence generated by PIAG is denoted by $\{x_k\}$. $\tilde x$ is said to be a stationary point of (\ref{model}) if $0\in\sum\nabla  f_i(\tilde x)+\partial h(\tilde x)$. The set of all stationary points of (\ref{model}) is denoted by $\mathcal{X}$. We say a sequence $\{x_k\}$ is R-linearly converges to $x^\ast$ if
$\limsup_{k\rightarrow+\infty}\|x_k-x^\ast\|^\frac{1}{k}<1$. 
%Each component function $f_i$ is assumed to have the $L_{f_i}$-Lipschitz continuous gradient ($f_i$ is $L_{f_i}$-smooth), i.e.,
%\begin{align}
%\|\nabla f_i(x)-\nabla f_i(y)\|\leq L_{f_i}\|x-y\|,~~\forall x,y\in\mathbb{R}^d.
%\end{align}
A function $f:\mathbb{R}^d\rightarrow \mathbb{R}$ is said to have a $L$-Lipschitz continuous gradient or to be $L$-smooth if
\begin{align}
\|\nabla f(x)-\nabla f(y)\|\leq L\|x-y\|,~~\forall x,y\in\mathbb{R}^d.
\end{align}
For the $L$-smooth function $f$, there always exists convex and gradient-Lipschitz continuous $f^{(j)},j=1,2$ such that $f=f^{(1)}-f^{(2)}$. As illustrated in \cite{peng}, one can choose $c>L$ and decompose $f$ in the following form:

\begin{align}
f= \underbrace{\left(f+\frac{c\|x\|^2}{2}\right)}_{f^{(1)}}-\underbrace{\frac{c\|x\|^2}{2}}_{f^{(2)}}.
\end{align}
We list assumptions involved in this paper as follows.
\begin{enumerate}
	\item[A0.] The objective function $F(x)$ in (\ref{model}) is lower bounded.
	\item[A1.] The decomposition $f_i=f_i^{(1)}-f_i^{(2)}$ exists for $i=1,\cdots,N$ such that $f_i^{(1)}$ is $L_i$-smooth and convex as well as $ f_i^{(2)}$ is $l_i$-smooth and convex. Denote $L=\sum_{i=1}^N L_i$ and $l = \sum_{i=1}^N l_i$. Also assume $L_i\geq l_i$ and thus $f_i$ is $L_i$-smooth.
	\item[A2.] The nonsmooth part $h:\mathbb{R}^d\rightarrow (-\infty,\infty]$ is proper, closed, convex and suqbdifferentiable everywhere in its effective domain, i.e., $\partial h(x)\neq\emptyset$ for all $ x\in\{y\in\mathbb{R}^d:h(y)<\infty\}$.
	\item[A3.] The time-varying delays $\tau_{k}^i$ are bounded; that is, there exists a nonnegative integer $\tau$ such that $\forall k\geq 1,i\in\{1,2,\cdots,N\}$, we have
	\begin{align}
	\tau_k^i\in\left\{0,1,\cdots,\tau \right\},
	\end{align}
	where $\tau$ is named as the \textit{delay parameter}.
\end{enumerate}	
The following two assumptions A4 and A5 are standard in the convergence analysis of several algorithms; see \cite{luo1992linear,luo1993convergence,beck2006linearly,TP1,TP3} and references therein. 
\begin{enumerate}
	\item[A4.] (\textbf{Proximal Error Bound Condition}) For any $\zeta\geq \inf_{x\in\mathbb{R}^d}F(x)$, there exist $\epsilon>0$ and $c_0>0$ such that
	\begin{align}
	\dist(x,\mathcal{X})\leq c_0\left\|\prox_{\frac{1}{L} h}\left(x-\frac{1}{L}\nabla f(x)\right)-x\right\|,
	\end{align}
	whenever $\|\prox_{\frac{1}{L}g}(x-\frac{1}{L}\nabla f(x))-x\|<\epsilon$ and $F(x)\leq\zeta$.
	\item[A5.] There exists $\delta>0$, such that $\|x-y\|\geq \delta$ whenever $x,y\in\mathcal{X},F(x)\neq F(y)$.
\end{enumerate}

For simplicity, we further denote 
$f=\sum_{i=1}^N f_i$, 
 $\bar l=\frac{l(\tau+1)}{2},\bar L=\frac{L(\tau+1)}{2}$ and  $\Delta_k=\sum_{j=k-\tau}^{k-1}\|x_{j+1}-x_j\|^2$. Without any loss of generality, let $x_{-k}=x_0$ for $k\geq 1$.
\section{Convergence Analysis}
First, we give a sufficient descent property of PIAG for nonconvex minimization.
\begin{lemma}\label{lemma1}With the assumptions A1-A3, the following statements for the problem (\ref{model}) hold:\\
\begin{enumerate}
	\item[(i)] For any $x\in\mathop{dom} F$, we have the descent lemma that
	\begin{align}
	F(x_{k+1})&\leq F(x)+\left(\bar l+\frac{1}{2\alpha}\right)\|x-x_k\|^2-\frac{1}{2\alpha}\|x-x_{k+1}\|^2\nonumber\\
	&+\left(\bar L-\frac{1}{2\alpha}\right)\|x_{k+1}-x_{k}\|^2+(\bar l+\bar L)\Delta_k,k\geq 0.
	\end{align}
	
	\item[(ii)]Consequently, we have the sufficient descent property that
	\begin{align}\label{lemma1_ineq4}
	F(x_{k+1})\leq F(x_k)+\left(\bar L-\frac{1}{\alpha}\right)\|x_{k+1}-x_k\|^2+(\bar l+\bar L)\Delta_k,k\geq 0.
	\end{align}
\end{enumerate}
\end{lemma}
\textit{Proof}. By the convexity of $f_i^{(1)}$ and the $l_{i}$-smoothness of $f_i^{(2)}$, we have
\begin{subequations}
\begin{align}
f_i^{(1)}(x)+\langle \nabla f^{(1)}_i(x),y-x\rangle&\leq f_i^{(1)}(y),\label{14a}\\
f_i^{(2)}(y)\leq f_i^{(2)}(x)+\langle \nabla f^{(2)}_i(x),&y-x\rangle+\frac{l_i}{2}\|y-x\|^2.\label{14b}
\end{align}
\end{subequations}
Adding up (\ref{14a}) and (\ref{14b}), using $f_i=f_i^{(1)}-f_i^{(2)}$, we obtain
\begin{align}\label{ineq1}
f_i(x)+\langle \nabla f_i(x),y-x\rangle\leq f_i(y)+\frac{l_i}{2}\|y-x\|^2.
\end{align}
Due to the $L_i$-smoothness of $f_i$ and the inequality (\ref{ineq1}), we have
\begin{align}\label{lemma1_ieqn1}
f_i(x_{k+1})&\leq f_i(x_{x-\tau_k^i})+\langle\nabla f_i(x_{k-\tau_k^i}),x_{k+1}-x_{k-\tau_k^i}\rangle +\frac{L_i}{2}\|x_{k+1}-x_{k-\tau^i_k}\|^2_2\nonumber\\
&\leq f_i(x)+\langle \nabla f_i(x_{k-\tau^i_k}),x_{k+1}-x\rangle+\frac{l_i}{2}\|x-x_{k-\tau^i_k}\|^2_2+\frac{L_i}{2}\|x_{k+1}-x_{k-\tau^i_k}\|^2_2.
\end{align}
Using the convexity of $\|\cdot\|^2$, we derive that
\begin{align}
\sum_{i=1}^N \frac{l_i}{2}\|x-x_{k-\tau_k^i}\|^2&= \sum_{i=1}^N \frac{l_i}{2}\|(x-x_k)+(x_k-x_{k-1})+\cdots+(x_{k-\tau_k^i+1}-x_{k-\tau_k^i})\|^2\nonumber\\
&\leq \sum_{i=1}^N\frac{l_i(\tau+1)}{2}\left(\|x-x_k\|^2+\sum_{j=k-\tau}^{k-1}\|x_{j+1}-x_j\|^2\right)\nonumber\\
&=\bar l\|x-x_k\|^2+\bar l\Delta_k\label{lineq}.
\end{align}
Similarly,
\begin{align}
\sum_{i=1}^N\frac{L_i}{2}\|x_{k+1}-x_{k-\tau_k^i}\|^2\leq \bar L\|x_{k+1}-x_k\|^2+\bar L\Delta_k.\label{Lineq}
\end{align}
With (\ref{lineq}) and (\ref{Lineq}), the sum of (\ref{lemma1_ieqn1}) from $i=1$ to $N$ becomes
\begin{align}\label{lemma1_ieqn3}
f(x_{k+1})\leq& f(x)+\langle g_k,x_{k+1}-x\rangle+\bar l\|x-x_k\|^2\nonumber\\
&+\bar L\|x_{k+1}-x_k\|^2+\bar l\Delta_k+ \bar L\Delta_{k}
\end{align}
where $g_k=\sum_{i=1}^N\nabla f(x_{k-\tau_k^i})$.
From the $\frac{1}{\alpha}-$strongly convexity of subproblem (\ref{intro_subprob}), we have
\begin{align}\label{lemma1_ieqn2}
\langle g_k,x_{k+1}-x\rangle\leq& h(x)-h(x_{k+1})+\frac{1}{2\alpha}\|x-x_k\|^2\nonumber\\
&-\frac{1}{2\alpha}\|x_{k+1}-x\|^2-\frac{1}{2\alpha}\|x_{k+1}-x_{k}\|^2.
\end{align}
Plugging (\ref{lemma1_ieqn2}) into (\ref{lemma1_ieqn3}), we obtain
\begin{align}
F(x_{k+1})
&\leq F(x)+\left(\bar l+\frac{1}{2\alpha}\right)\|x-x_k\|^2-\frac{1}{2\alpha}\|x-x_{k+1}\|^2\nonumber\\
&+\left(\bar L-\frac{1}{2\alpha}\right)\|x_{k+1}-x_{k}\|^2+(\bar l+\bar L)\Delta_k.
\end{align}
Then the statement (i) holds. The statement (ii) follows from  statement (i) by setting $x=x_k$.
\qed

Through the sufficient descent property of nonconvex PIAG, we give the following lemma to illustrate that for a fixed positive integer $M$, the sequence $\{x_k\}$ satisfies
 \begin{align}
 \lim_{k\rightarrow\infty} \|x_{k+M}-x_k\|=\lim_{k\rightarrow\infty}\sum_{j=k}^{j=k+M}\|x_{j+1}-x_j\|\rightarrow0.
 \end{align}

%Inspired by monotonicity of the project gradient operator(see lemma 25 in \cite{wang2014iteration}), which is an special case of the PG method by regarding $h(x)$ as an indicator function, we claim that  the PG method is also monotonically increasing respect to the stepsize. Since 
%To make it explicit, we give the following lemma.

\begin{lemma}\label{lemma2} Assume A0-A3 hold. If stepsize $\alpha<\frac{1}{\bar L+\tau (\bar l+\bar L)}$, then the following statements hold:
	\begin{enumerate}
		\item[(i)] $F(x_k)$ is bounded;
		\item[(ii)] $\sum_{k=0}^\infty\|x_{k+1}-x_k\|^2< +\infty$.
	\end{enumerate}
\end{lemma}
\textit{Proof.} 
From Lemma \ref{lemma1}(ii), for arbitrary positive integers $k_1<k_2$,  summing up (\ref{lemma1_ineq4}) from $k=k_1$ to $k_2-1$ yields
\begin{align}\label{lemma2_iqn}
F(x_{k_2})&\leq F(x_{k_1})+\left(\bar L-\frac{1}{\alpha}\right)\sum_{k=k_1}^{k_2-1}\|x_{k+1}-x_k\|^2+(\bar l+\bar L)\sum_{k=k_1}^{k_2-1}\Delta_k\nonumber\\
&\leq F(x_{k_1})+\left(\bar L-\frac{1}{\alpha}\right)\sum_{k=k_1}^{k_2-1}\|x_{k+1}-x_k\|^2+\tau(\bar l+\bar L)\sum_{k=k_1-\tau}^{k_2-1}\|x_{k+1}-x_k\|^2\nonumber\\
&\leq F(x_{k_1})+\left(\bar L+\tau(\bar l+\bar L)\right)\sum_{k=k_1-\tau}^{k_2-1}\|x_{k+1}-x_k\|^2-\frac{1}{\alpha}\sum_{k=k_1}^{k_2-1}\|x_{k+1}-x_k\|^2.
\end{align}
Setting $k_1=0$ and $k_2=K+1$, we obtain
\begin{align}\label{lemma2_ieqn}
F(x_{K+1})\leq F(x_{0})+\left(\bar L-\frac{1}{\alpha}+\tau(\bar l+\bar L)\right)\sum_{k=0}^{K}\|x_{k+1}-x_k\|^2,
\end{align}
which indicates that $F(x_k)$ is bounded from above if $\alpha<\frac{1}{\bar L+\tau (\bar l+\bar L)}$.  With A0 that $\inf F> -\infty$ holds, (\ref{lemma2_ieqn}) implies

\begin{align}
\sum_{k=0}^K\|x_{k+1}-x_k\|^2\leq\frac{F(x_0)-F(x_{K+1})}{\frac{1}{\alpha}-\tau(\bar l+\bar L)-\bar L}.
\end{align}
The inequality holds as $K\rightarrow \infty$. Thus statement (ii) is proved.
\qed

\begin{lemma}\label{lemma3}
	Assume that A0-A3 hold and $\alpha<\frac{1}{\bar L+\tau (\bar l+\bar L)}$. Then, any accumulation point of $\{x_k\}$ is a stationary point of $F$.
\end{lemma}
\textit{Proof.} Let $\bar x$ be an accumulation point. Then there exists a subsequence $\{x_{k_i}\}$ such that $\lim_{i\rightarrow\infty}x_{k_i}=\bar x$. Using the first-order optimality condition of subproblem (\ref{intro_subprob}), we have
\begin{align}\label{belong}
-\frac{1}{\alpha}({x_{k_i+1}}-x_{k_i})\in\sum_{j=1}^N \nabla f_j(x_{k_i-\tau_{k_i}^j})+\partial h(x_{k_i+1}).
\end{align}
Invoking Lemma 2(ii), for an arbitrary fixed integer $I\in\{0,1,\cdots,\tau\}$, we have
\begin{align}
\lim_{i\rightarrow+\infty} x_{k_i-I}=\lim_{i\rightarrow+\infty} x_{k_i+1}=\bar x,
\end{align}
which implies $\lim_{i\rightarrow+\infty} x_{k_i-\tau_{k_i}^j}=\bar x$ since $0\leq\tau_{k_i}^j\leq\tau$ for any $j\in\{1,2,\cdots,N\}$.
Consequently, due to continuity of $\nabla f$ and closedness of $\partial h$, (\ref{belong}) implies $0\in\nabla f(\bar x)+\partial h(\bar x)$.
\qed
\begin{lemma}\label{lemma4}Assume that A0-A3 hold and $\alpha<\frac{1}{\bar L+\tau (\bar l+\bar L)}$. Let $\Omega$ be the set of accumulation points of the sequence $\{x_k\}$ generated by PIAG. Then $\zeta=lim_{k\rightarrow \infty} F(x_k)$ exists and $F\equiv\zeta$ on $\Omega$.
\end{lemma}
\textit{Proof.} The fact that $F(x_k)$ is bounded has been shown in Lemma \ref{lemma2}(i). Thus if the limit of $F(x_k)$ does not exist, then there are two subsequences $\{x_{s^1_i}\}$ and $\{x_{s^2_i}\}$ of $\{x_k\}$ such that $F(x_{s_i^{1}})\rightarrow F_1$ and $F(x_{s_i^{2}})\rightarrow F_2$. Without loss of generality, suppose $F_1>F_2$. 

First, due to Lemma \ref{lemma2}(ii), there exists a sufficiently large positive integer $K$ such that
\begin{align}
\left(\bar L+\tau(\bar l+ \bar L)\right)\sum_{j=K-\tau}^\infty\left\|x_{j+1}-x_{j}\right\|^2< \frac{F_1-F_2}{3}.
\end{align}

Second, from (\ref{lemma2_iqn}) we can find two sufficiently large subscript indexes $S_1\in\{s_i^{1}\}$ and $S_2\in\{{s_{i}^{2}}\}$ such that $S_1> S_2>K$ and satisfy
\begin{align}
F(x_{S_1})-F(x_{S_2})&\leq \left(\bar L+\tau(\bar l+ \bar L)\right)\sum_{j=S_2-\tau}^{S_1-1}\|x_{j+1}-x_j\|^2\nonumber\\
&< \frac{F_1-F_2}{3}\label{ieqn1},\\
F_1&-F(x_{S_1})< \frac{F_1-F_2}{3}\label{ieqn2},\\
F&(x_{S_2})-F_2< \frac{F_1-F_2}{3}\label{ieqn3}.
\end{align}
The sum of (\ref{ieqn1}),(\ref{ieqn2}) and (\ref{ieqn3}) derives the contradiction $F_1-F_2< F_1-F_2$. Thus $\lim F(x_k)$ must exist. Denote the limit by $\zeta$.

If $\Omega=\emptyset$, the result is trivially true. Otherwise, $\forall \hat x\in\Omega$, suppose a subsequence $x_{k_i}\rightarrow\hat x$. Due to the lower semi-continuity of $F$, we have
\begin{align}\label{lemma4_ieqn}
F(\hat x)\leq \lim\inf F(x_{k_i})=\zeta.
\end{align}

On the other hand, since $x_{k_i+1}$ is the minimizer of 
\begin{align}
h(x)+\langle g_{k_i},x-x_{k_i}\rangle+\frac{1}{2\alpha}\|x-x_{k_i}\|^2,
\end{align}
we have
\begin{align}
&f(x_{k_i+1})+h(x_{k_i+1})+\langle g_{k_i},x_{k_i+1}-\hat x\rangle+\frac{1}{2\alpha}\| x_{k_i+1}-x_{k_i}\|^2\\
\leq& f(x_{k_i+1})+h(\hat x)+\frac{1}{2\alpha}\|\hat x-x_{k_i}\|^2.
\end{align}
Letting $i\rightarrow \infty$, we obtain
\begin{align}
\lim\sup F(x_{k_i})\leq F(\hat x).
\end{align}
Along with (\ref{lemma4_ieqn}), the equality $F(\hat x)\equiv\zeta$ holds for all $\hat x\in\Omega$.
\qed
Furthermore, if we assume $F$ is level bounded, since we already know $F(x_k)$ is upper bounded from Lemma \ref{lemma2}, then the sequence $\{x_k\}$ is also bounded which implies$\Omega\neq\emptyset$ in l  emma \ref{lemma3} implies that .

The stepsize $\alpha$ is required to be small in previous lemmas but is undetermined for now. We might require a sufficiently small $\alpha$ in PIAG to guarantee convergence. The A4 for a fixed stepsize $\alpha=\frac{1}{L}$ seems inadequate for later proof. Therefore, we need a variant of A4 with stepsizes smaller than $\frac{1}{L}$. To make the fact explicit, we display the result in the following two lemmas.

\begin{lemma}\cite[lemma 2]{RN83}\label{monoLemma} Suppose that function $h:\mathbb{R}^d\rightarrow \mathbb{R}$ satisfies A2 and $f:\mathbb{R}^d\rightarrow \mathbb{R}$ is differentiable on $\mathbb{R}^d$. Then, $\forall x\in \dom h$ and real numbers $t\geq t'>0$, we have
	\begin{align}
	 \frac{1}{t}\|\prox_{t h}(x-t\nabla f(x))-x\|\leq \frac{1}{t'}\|\prox_{t h}(x-t'\nabla f(x))-x\|.
	\end{align}
	 
\end{lemma}

\begin{lemma}\label{equiv} If A4 holds with $\epsilon>0,c_0>0$ and $\zeta\geq \inf_{x\in\mathbb{R}^d}F(x)$, then
	\begin{align}\label{proxeb}
	\dist(x,\mathcal{X})\leq \frac{c_0}{\alpha L}\left\|\prox_{\alpha h}\left(x-\alpha\nabla f(x)\right)-x\right\|,
	\end{align}
	whenever $0<\alpha\leq\frac{1}{L}$, $\|\prox_{\alpha h}(x-\alpha\nabla f(x))-x\|<\alpha L\epsilon$ and $F(x)\leq\zeta$.
\end{lemma}
\textit{Proof.}
For $\alpha\in(0,\frac{1}{L}]$, if $x$ satisfies $\|\prox_{\alpha h}(x-\alpha\nabla f(x))-x\|<\alpha L\epsilon$ and $F(x)\leq\zeta$, invoking Lemma \ref{monoLemma}, then we have
\begin{align}
\|\prox_{\frac{1}{L} h}(x-\frac{1}{L}\nabla f(x))-x\|\leq\frac{1}{\alpha L}\|\prox_{\alpha h}(x-\alpha\nabla f(x))-x\|<\epsilon.
\end{align}
Thus A4 gives
\begin{align}
\dist(x,\mathcal{X})&\leq c_0\left\|\prox_{\frac{1}{L} h}\left(x-\frac{1}{L}\nabla f(x)\right)-x\right\|\\
&\leq \frac{c_0}{\alpha L}\left\|\prox_{\alpha h}\left(x-\alpha\nabla f(x)\right)-x\right\|,
\end{align}
which is just (\ref{proxeb}).
\qed
Before proposing the final results, we need the following key lemma for revealing the linear convergence.
\begin{lemma}\cite{PIAG2016}\label{lemma_bound}
	Assume that the non-negative sequences $\{V_k\}$ and $\{\omega_k\}$ satisfy the following inequality:
	\begin{align}
	V_{k+1}\leq a V_k-b\omega_k+c\sum_{j=k-k_0}^k\omega_j,
	\end{align}
	for some real numbers $a\in(0,1)$ and $b,c\geq 0$, and some positive integer $k_0$. Also Assume that $\omega_k=0$ for $k<0$, and that the following holds:
	\begin{align}\label{condtion_bound}
	\frac{c}{1-a}\frac{1-a^{k_0+1}}{a^{k_0}}\leq b.
	\end{align}
	Then, $V_k\leq a^kV_0$ for all $k\geq 0$.
\end{lemma}
\begin{theorem}
For arbitrary sequence $\{x_k\}$ generated by PIAG, with the assumption A0-A5, for all sufficiently small stepsize $\alpha\geq 0$, the following statements hold:
\begin{enumerate}
	\item[(i)]$\dist(x_k,\mathcal{X})\rightarrow 0$;
	\item[(ii)]$F(x_k)$ is R-linearly convergent.
\end{enumerate}
\end{theorem}
\textit{Proof}. We derive that
\begin{align}
&\|\prox_{\alpha h}(x_k-\alpha\nabla f(x_k))-x_k\|\nonumber\\
\leq&\|\prox_{\alpha h}(x_k-\alpha \nabla f(x_k))-\prox_{\alpha h}(x_k-\alpha g_k)\|+\|\prox_{\alpha h}(x_k-\alpha g_k)-x_k\|\nonumber\\
\leq&\alpha\sum\|\nabla f_i(x_k)-\nabla f_i(x_{k-\tau_k^i})\|+\|x_{k+1}-x_k\|\nonumber\\
\leq&\alpha\sum L_i\|x_k-x_{k-\tau_{k}^i}\|+\|x_{k+1}-x_k\|\nonumber\\
\leq&\sqrt{(\alpha L+1)\left(\alpha\sum_{i=1}^N L_i\|x_k-x_{k-\tau_k^i}\|^2+\|x_{k+1}-x_k\|^2\right)}\nonumber\\
\leq&\sqrt{(aL+1)\left(\alpha\tau L\sum_{j=k-\tau}^{k-1}\|x_{j+1}-x_j\|^2+\|x_{k+1}-x_{k}\|^2\right)}\label{proxleq}
\end{align}
In addition, we require 
\begin{align}
\alpha<\min\left(\frac{1}{\bar L+\tau (\bar l+\bar L)},\frac{1}{L}\right).
\end{align} 
With Lemma \ref{lemma2}(ii), the inequality (\ref{proxleq}) implies
\begin{align}\label{proxto0}
\|\prox_{\alpha h}(x_k-\alpha\nabla f(x_k))-x_k\|\rightarrow 0.
\end{align}
In addition to (\ref{proxto0}), Lemma \ref{lemma4} implies that $F(x_k)$ is bounded, so we can conclude there exists a large enough positive number $K_\alpha$ relevant to $\alpha$ such that the inequality (\ref{proxeb}) could apply whenever $k>K_\alpha$. Therefore for each $k>K_\alpha$, using Assumption A4 and inequality (\ref{proxleq}) we have 
\begin{align*}
dist(x_k,\mathcal{X})^2
&\leq \left(\frac{c_0}{\alpha L}\right)^2(aL+1)\left(\alpha\tau L\sum_{j=k-\tau}^{k-1}\|x_{j+1}-x_j\|^2+\|x_{k+1}-x_{k}\|^2\right),
\end{align*}
which readily leads to $\dist(x_k,\mathcal{X})\rightarrow 0$. Thus the statement (i) is proved.

Setting $x=\bar x_k$ in (\ref{lemma1_ieqn1}), where $\bar x_k$ represents a projection of $x_k$ onto $\mathcal{X}$ then we obtain
\begin{align}
F(x_{k+1})&\leq F(\bar x_k)+\left(\frac{l(\tau+1)}{2}+\frac{1}{2\alpha}\right)dist(x_k,\mathcal{X})^2+\left(\frac{L(\tau+1)}{2}-\frac{1}{2\alpha}\right)\|x_{k+1}-x_{k}\|^2\nonumber\\
&+\frac{(l+L)(\tau+1)}{2}\sum_{j=k-\tau}^{k-1}\|x_{j+1}-x_j\|^2\nonumber\\
&\leq F(\bar x_{k}) +\left(\frac{L(\tau+1)}{2}-\frac{1}{2\alpha}+\left(\frac{l(\tau+1)}{2}+\frac{1}{2\alpha}\right) \left(\frac{c_0}{\alpha L}\right)^2(\alpha L+1)\right)\|x_{k+1}-x_{k}\|^2\nonumber\\
&+\left[\frac{(l+L)(\tau+1)}{2}+\left(\frac{l(\tau+1)}{2}+\frac{1}{2\alpha}\right) \left(\frac{c_0}{\alpha L}\right)^2(\alpha L+1)\alpha\tau L\right]\sum_{j=k-\tau}^{k-1}\|x_{j+1}-x_j\|^2\label{bigone}
\end{align}
Since $\|x_k-x_{k+1}\|\rightarrow 0$ and $\dist(x_k,\mathcal{X})\rightarrow 0$, nothing that
\begin{align}
\|\bar x_k-\bar x_{k+1}\|\leq \dist(x_k,\mathcal{X})+\|x_k-x_{k+1}\|+\dist(x_{k+1},\mathcal{X}),
\end{align}
we have $\|\bar x_k-\bar x_{k+1}\|\rightarrow 0$. With A5, $F(\bar x_k)\equiv \zeta$ holds for some constant $\zeta$ for all sufficiently large $k$. Without loss of generality, we suppose $F(\bar x_{k})=\zeta$ for $k\geq K_\alpha$. Relax $\alpha$ to $\frac{1}{L}$ and rewrite (\ref{lemma1_ineq4}) to (\ref{relate}), $\alpha^2\times$(\ref{bigone}) to (\ref{relate2}) in a simplified way as follow:
\begin{align}
&F(x_{k+1})\leq F(x_k)+\left(C_1-\frac{1}{\alpha}\right)\|x_{k+1}-x_k\|^2+C_2\sum_{j=k-\tau}^{k-1}\|x_{j+1}-x_j\|^2\label{relate},\\
&\alpha^2 F(x_{k+1})\leq \alpha^2\zeta+(C_3+\frac{c_0^2}{L^2}\frac{1}{2\alpha})\|x_{k+1}-x_k\|^2+C_4\sum_{j=k-\tau}^{k-1}\|x_{j+1}-x_j\|^2\label{relate2},
\end{align}
where constants $C_i$ are independent with $\alpha$. 
Denote $H(x_k):=F(x_k)-\zeta$. Therefore, $(\ref{relate})+\frac{L^2}{c_0^2}(\ref{relate2})$ leads to
\begin{align}
\left(1+\frac{L^2}{c_0^2}\alpha^2\right)H(x_{k+1})\leq H(x_k)+\left(C_5-\frac{1}{2\alpha}\right)\|x_{k+1}-x_k\|^2+C_6\sum_{j=k-\tau}^{k-1}\|x_{j+1}-x_j\|^2.\nonumber
\end{align}
Actually, from the inequality above, one can directly conclude that for all sufficiently small $\alpha$, Lemma \ref{lemma_bound} could always be employed to obtain the linear convergence of $H(x_k)\rightarrow 0$. The remaining trivial piece is to give an explicit range of $\alpha$.
All constants are listed as follows (a verifying Walfram Mathematica script is available online https://www.deepinfar.cn/piag):
\begin{align}
C_1&=\frac{L(\tau+1)}{2},C_2=\frac{(l+L)(\tau+1)}{2},\nonumber\\
C_3&=\frac{c_0^2 (2 l (\tau+1)+L)+L \tau}{2 L^2},\nonumber\\
C_4&=\frac{(l + L) (1 + \tau) + 2 \tau (l + L + l \tau) 
c_0^2}{2 L^2},\nonumber\\
C_5&=l + L + l \tau + \frac{L \tau}{2} + \frac{L\tau}{2c_0^2},\nonumber\\
C_6&=\frac{1}{2} \left(\frac{(\tau +1) (l+L)}{c_0^2}+2 l \tau ^2+3 l \tau +l+3 L \tau +L\right).\nonumber
\end{align}
With consistent notations in Lemma \ref{lemma_bound}, let $V_k=\left(1+\frac{L^2}{c_0^2}\alpha^2\right)H(x_k)$, $a=\left(1+\frac{L^2}{c_0^2}\alpha^2\right)^{-1}$, $b=\frac{1}{2\alpha}-C_5$,  $k_0=\tau$ and $c=C_6$. Note $\alpha\leq\frac{1}{L}$, then $a\geq\left(1+\frac{1}{c_0^2}\right)^{-1}$ and the left side of (\ref{condtion_bound}) is bounded by
\begin{align}
\frac{c}{1-a}\frac{1-a^{k_0+1}}{a^{k_0}}=c\sum_{j=0}^{k_0}a^{-j}\leq  C_6(1+\tau(1+\frac{1}{c_0^2})^\tau)\triangleq C_7,
\end{align}
which implies that when $\alpha\leq\frac{1}{2C_5+2C_7}$, we have
\begin{align}
\frac{c}{1-a}\frac{1-a^{k_0+1}}{a^{k_0}}\leq C_7\leq \frac{1}{2\alpha}-C_5=b.
\end{align}
 The inequality (\ref{condtion_bound}) in Lemma \ref{lemma_bound} holds. Thus let
 \begin{align}
 \alpha=C_8=\min\left(\frac{1}{\bar L+\tau (\bar l+\bar L)},\frac{1}{2C_5+2C_7},\frac{1}{L}\right)
 \end{align}
and we have
\begin{align}
F(x_k)-\zeta\leq\left(1+\frac{L^2}{c_0^2}\frac{1}{C_8^2}\right)^{K_\alpha-k}\left(F(x_{K_\alpha})-\zeta\right),~~k\geq K_\alpha
\end{align}
\qed
Moreover, we claim that the path of $\{x_k\}$ has finite length
and $\{x_k\}$ R-linearly converges to some stationary point $\tilde x$. First, we prove the following lemma.
\begin{lemma}\label{lemmaa}Let $\{a_k\},\{b_k\}$ be positive sequences where $b_k=b_0q^k$ for a real number $q\in(0,1)$. If the inequality
	\begin{align}
	a_k\leq b_k+\frac{c}{\tau}\left(a_{k-1}+a_{k-2}+\cdots+a_{k-\tau}\right)
	\end{align}
	holds for $k\geq \tau$ where $\tau$ is a given positive integer and $0<c<1$, then $\{a_k\}$ is R-linearly convergent to zero.
\end{lemma}
\textit{Proof.} We consider a characteristic polynomial
\begin{align}
P(x)=x^\tau-\frac{c}{\tau}x^{\tau-1}-\frac{c}{\tau}x^{\tau-2}-\cdots-\frac{c}{\tau}.
\end{align}
Since $P(c)\leq 0$ and $P(1)=1-c>0$, letting a root of $P(x)$ in $[c,1)$ be denoted by $p$,
then we have the following inequalities,
\begin{align*}
a_k&\leq b_k+\frac{c}{\tau}\left(a_{k-1}+a_{k-2}+\cdots+a_{k-\tau}\right),\\
pa_{k-1}&\leq pb_{k-1}+\frac{c}{\tau}\left(pa_{k-2}+pa_{k-3}+\cdots+pa_{k-\tau-1}\right),\\
p^2a_{k-2}&\leq p^2b_{k-2}+\frac{c}{\tau}\left(p^2a_{k-3}+p^2a_{k-4}+\cdots+p^2a_{k-\tau-2}\right),\\
&\cdots\\
p^{k-\tau}a_{\tau}&\leq p^{k-\tau}b_{\tau}+\frac{c}{\tau}\left(p^{k-\tau}a_{\tau-1}+p^{k-\tau}a_{\tau-2}+\cdots+p^{k-\tau}a_{0}\right).
\end{align*}
Adding them up, with proper relaxing, we have
\begin{align}
a_k&\leq \sum_{i=0}^k p^{k-i}b_i+cp^{k-2\tau+1}\sum_{i=0}^{\tau-1}a_i\\
&=b_0p^k\sum_{i=0}^k\left(\frac{q}{p}\right)^{i}+cp^{k-2\tau+1}\sum_{i=0}^{\tau-1}a_i.
\end{align}
Perform limit superior on both sides, we have
\begin{align}
\mathop{\lim\sup}_{k\rightarrow +\infty}{|a_k|^{\frac{1}{k}}}\leq\max(p,q)<1.
\end{align}
Thus $\{a_k\}$ is R-linearly convergent to $0$.\qed

\begin{theorem}
	Suppose conditions of Theorem 1 are satisfied. Then we have
	\begin{enumerate}
		\item[(i)] $\sum_{i=0}^\infty \|x_{i+1}-x_i\|<+\infty$,
		\item[(ii)] $\{x_k\}$ is R-linearly convergent.
	\end{enumerate}
\end{theorem}
\textit{Proof.}
Rewrite the inequality (\ref{relate}) as
\begin{align}
\|x_{k+1}-x_k\|^2\leq \frac{C_2}{\frac{1}{\alpha}-C_1}\left(F(x_k)-F(x_{k+1})\right)+\frac{C_2}{\frac{1}{\alpha}-C_1}\sum_{j=k-\tau}^{k-1}\|x_{j+1}-x_{j}\|^2.
\end{align}
Since the first term of the right sides is proved to be R-linearly convergent and the coefficient of the second term satisfies
\begin{align}
\frac{C_2}{\frac{1}{\alpha}-C_1}=\frac{C_2}{2C_5+2C_7-C_1}<\frac{C_2}{2C_5+2\tau C_6-C_1}<\frac{1}{\tau},
\end{align}
thus Lemma \ref{lemmaa} implies that 
$\|x_{k+1}-x_k\|\leq r_0\cdot r^k$ for some $0<r<1,r_0>0$, which illustrates the statement (i). Consequently $\{x_k\}$ is a Cauchy sequence and then it converges to a point $\tilde x\in\mathcal{X}$. Finally we have
\begin{align}
\|x_{k}-\tilde x\|\leq\sum_{j=k}^\infty \|x_{j}-x_{j+1}\|\leq\frac{r_0\cdot r^k}{1-r},
\end{align}
which implies that $x_k$  R-linearly converges to $\tilde x$.
\qed

\section{Conclusion}
In this paper, we analyze the convergence of PIAG for nonconvex minimization. First of all, we give the sufficient descent property of PIAG in nonconvex cases. Under the proximal error bound condition, we prove that the generated sequence $\{x_k\}$ is convergent to the stationary point set. Then, we show $\{F(x_k)\}$ is R-linearly convergent and that $\{x_k\}$ R-linearly converges to a stationary point when the stepsize $\alpha$ is under some positive constant. Finally, we note that even with the delay parameter $\tau$ vanishing, our theoretical convergence rate is far from being tight, which deserves further study.

%\begin{acknowledgements}
%\end{acknowledgements}


%\begin{thebibliography}{}
% etc
%\end{thebibliography}

\bibliographystyle{plain}
\bibliography{ref}

\end{document}